\documentclass[12]{article}

\usepackage[cp1251]{inputenc}
\usepackage[english]{babel}
\usepackage{latexsym}
\usepackage{amscd}
\usepackage{amssymb}
\usepackage{amsmath}
\usepackage{amsthm}

\newtheorem{lm}{Lemma}
\newtheorem{theo}{Theorem}
\newtheorem{rem}{Remark}

\title{ On the Structure of the Bochner--Martinelli \\
Residue Currents
} 

\author
{ Irina Antipova\thanks{ The research was carried out in Siberian Federal University within the framework of the project, funded by the grant of the Russian Federation Government to support scientific research under supervision of a leading scientist, no. 14.Y26.31.0006. The author was also supported by Russian Foundation for Basic Research ( 13-01-12417 ofi-m2,  14-01-00544 A) .} }
\date{}

\begin{document}
\large
{
\maketitle

\begin{abstract}
We study residue currents of the Bochner--Martinelli type using their relationship  with Mellin transforms of  residue integrals.
We present the structure formula for residue currents associated with monomial mappings: they admit representations as sums of products of   residue currents, principal value currents and  hypergeometric functions.

\bigskip

{\bf Keywords}: residue current, Bochner--Martinelli form, Mellin transform.
\end{abstract}

\section{Introduction}

Let us consider a residue current of the Bochner--Martinelli type associated with a holomorphic mapping
\begin{equation}\label{f0}
f=(f_1, \dots , f_p):V\to {\mathbb C}^p,
\end{equation}
where $V$ is some domain in $ {\mathbb C}^n$ and $p\in\{1,\dots , n\}$. Such a current was defined by M. Passare, A. Tsikh and A. Yger in \cite{PTY} as a limit
\begin{equation}\label{f2}
T_{f}(\varphi)=\lim\limits_{\eta \to 0} c_p \int\limits_{\left\|f\right\|^2=\eta}^{}\Omega (f)\wedge \varphi,
\end{equation}
where the kernel $\Omega (f)$ is given by
$$
\Omega (f)= \frac{1}{\left\|f\right\|^{2p}} 
\sum\limits_{k=1}^{p} (-1)^{k-1} \overline{f_k}  \bigwedge_{l\ne k}^{} \overline{df_l}, 
$$
the test form  $\varphi$ with bidegree $(n,n-p)$  is smooth and compactly supported in $V$ ($\varphi \in {\mathcal D}^{n,n-p} (V)$) and 
$c_p=\frac{(-1)^{p(p-1)/2}(p-1)!}{(2\pi i)^p}$ is a constant.
It may be viewed as a limit of the mean value of the residue integral
\begin{equation}\label{f1}
I(\varepsilon)=I(\varepsilon, \varphi) =\frac{1}{(2\pi i)^{p}}\int\limits_{{\mathcal T}_{\varepsilon}(f)} 
\frac{\varphi}{f_1 \cdot\ldots \cdot f_p}
\end{equation}
over the simplex
$$
\Sigma_{p}(\eta)=\left\{\varepsilon\in {\mathbb R}^{p}_{+}: \varepsilon_1+\dots + \varepsilon_p=\eta\right\}.
$$
  The integration set in (\ref{f1}) is the real analytic chain (tube) 
	$$
	{\mathcal T}_{\varepsilon}(f)=\{|f_1|^{2}=\varepsilon_{1},\dots , |f_p|^{2}=\varepsilon_{p}\},
	$$
	which is oriented as the skeleton of the corresponding polyhedron.  Theorem 1.1 in \cite{PTY} states that  limit (\ref{f2})  always exists and defines a $(0, p)$--current. In the complete intersection case, when the set  $f^{-1}(0)$ has codimension  $p$, the current  $T_f$ coincides with the Coleff--Herrera current (see \cite{PTY, CH}) and with currents considered in papers  \cite{P}, \cite{PT95}, \cite{BGY}. The advantage of this approach is caused by the fact that the integral in (\ref{f2}) depends on one parameter $\eta$ instead of $p$ parameters  $\varepsilon_1,\dots ,\varepsilon_p$ in (\ref{f1}).

The concept of the Bochner-Martinelli residue currents has been considerably developed by M. Andersson, E. Wulcan and other authors; see, for example  \cite{A, W07}, for more details. The structural formula for the Bochner--Martinelli residue current $R^{z^{A}}$ associated with the monomial ideal  $(z^{A})$ of zero dimension, $A\subset {\mathbb Z}_{+}^{n}$, was studied in \cite{W07}. In particular, a description of the  annihilator ideal ${\mbox{\rm Ann}}R^{z^{A}}$ in terms of the associated Newton diagram was given. 
 
Our goal is to study the structure of the Bochner--Martinelli residue current in the normal crossing case, when  $f_j$ in (\ref{f0}) are monomials. Theorem 1 claims that the  current associated with the system of $p$ monomials in $n$ variables may be represented  as the sum of  products of residue currents, principal value currents and   hypergeometric functions. Essentially, we are interested in the case when $1<p<n$. For the case $p=1$ the problem was studied  by P. Dolbeault in \cite{D}. The Bochner--Martinelli current for a monomial ideal, generated by $n$ monomials (i.e. in the case $p=n$), 
was calculated in \cite{PTY}. In the present paper we  treat ideas proposed in \cite{PTY} in general monomial case. We focus on the demonstration of the Mellin transform technique (see \cite{A07}) and residue theory methods that have been proven in the theory of residue currents (see \cite{PT95, TY04}).

Let us remind that for a holomorphic function $f:V\to {\mathbb C}$ in some domain $V\subset {\mathbb C}^n$ the principal value of $1/f$, determined by the limit
$$
\left\langle \left[1/f \right],\psi\right\rangle=\lim\limits_{\varepsilon\to 0}^{}\int\limits_{\{|f|^2>\varepsilon\}}^{}\frac{\psi}{f},
$$
defines a current for $\psi\in {\mathcal D}^{n,n}(V)$ (see \cite{CH}). By Stokes' formula the $\bar{\partial}$--action of this current is equal to
$$
\langle \bar\partial [1/f],\varphi\rangle=-\lim\limits_{\varepsilon\to 0}^{}\int\limits_{\{|f|^2>\varepsilon\}}^{}\frac{\bar\partial\varphi}{f}=\lim\limits_{\varepsilon\to 0}^{}\int\limits_{\{|f|^2=\varepsilon\}}^{}\frac{\varphi}{f},\,\, {\varphi}\in {\mathcal D}^{n,n-1} (V).
$$
In particular, when $n=1$ and $f=z^k$ is a monomial, one has
$$
\left\langle \bar\partial \left[1/z^k\right],\varphi\right\rangle = \frac{2\pi i}{(k-1)!}\frac{\partial^{k-1}}{\partial z^{k-1}}\tilde{\varphi} (0),
$$ 
where $\varphi = \tilde{\varphi}(z) dz$.

\section{Role of the Mellin Transform}

Rewrite the integral in (\ref{f2}) in the form
\begin{equation*}
J(\eta)=J(\eta , \varphi)=\frac{c_p}{\eta^p}\int\limits_{\left\|f\right\|^2=\eta}^{}\left(\sum\limits_{k=1}^{p}(-1)^{k-1}\overline{{f}_k} \bigwedge_{l\ne k}^{} \overline{df_l}\right)\wedge \varphi ,
\end{equation*}
where $\varphi$ is a test form from the space ${\mathcal D}^{n,n-p} (V)$.
Its Mellin transform
\begin{equation}\label{f1.1}
M[J(\eta)](\varphi , s):=s\int\limits_{0}^{+\infty} J(\eta) \eta ^{s-1}d\eta
\end{equation}
is a holomorphic function in the half-plane $\Re s> p$. According to \cite{PTY} it is holomorphic across the imaginary axis and  can be meromorphically continued to the whole complex plane. Moreover, the map
\begin{equation*}
{\mathcal D}^{n,n-p}\ni\varphi \mapsto M[J(\eta)](\varphi , 0)
\end{equation*}
gives the value of the current $T_{f}$ on $V$. This current  can be realized also as the following limit of
solid volume integrals
\begin{equation}\label{f2.1}
T_{f}(\varphi) =\lim\limits_{\tau \to 0+} p c_p\int\limits_{V}^{}\frac{\tau \bigwedge_{k=1}^{p} \overline{{df}_{k}}\wedge \varphi}{\left(\left\|f\right\|^2+\tau\right)^{p+1}}. 
\end{equation} 
One can say that (\ref{f2.1}) is a regularization for $T_f$, since the integrand here is a regular differential form for any $\tau >0$.
Using representation (\ref{f2.1}), let us display the role of the Mellin transform of residue integral (\ref{f1}) in realization of the current $T_{f}(\varphi)$. To outline the idea we take $p=1$. Applying the  Mellin transform inversion formula for  the function
$$
\frac{\tau}{(|f|^2+\tau)^2}=:F(\tau),
$$
we obtain the following Mellin--Barnes representation
\begin{equation}\label{f2.2}
F(\tau)=\frac{1}{2\pi i}\int\limits_{\gamma-i\infty}^{\gamma+i\infty} \Gamma (s+1) \Gamma (1-s) \left|f\right|^{2(s-1)}\tau^{-s}ds,\,\, 0<\gamma <1.
\end{equation}
After substitution of the expression (\ref{f2.2}) in to integral (\ref{f2.1}) we get by Fubini's theorem the following representation:
\begin{equation}\label{f2.4}
T_{f}(\varphi) =\lim\limits_{\tau\to 0+}\frac{1}{2\pi i}\int\limits_{\gamma-i\infty}^{\gamma+i\infty}\Gamma (s+1)\Gamma (1-s)\Gamma_{f}(s,\varphi) \tau^{-s} ds,\,\, 0<\gamma<1,
\end{equation}
where  
\begin{equation*}\label{f2.3}
\Gamma_{f}(s,\varphi)=\frac{1}{2\pi i}\int\limits_{V}^{} |f|^{2(s-1)}\overline{df}\wedge\varphi .
\end{equation*}
The $\Gamma_{f}(s,\varphi) $ is equal to the Mellin transform 
$$
\int\limits_{0}^{+\infty}I(\varepsilon , \varphi) \varepsilon^{s-1} d\varepsilon
$$
of the residue function
$$
I(\varepsilon, \varphi) =\frac{1}{2\pi i}\int\limits_{|f|^2={\varepsilon}} 
\frac{\varphi}{f}.
$$

A multidimensional analog of formula (\ref{f2.4}) (see below (\ref{f5})) was obtained in \cite{VY} in the complete intersection case. In fact, this hypothesis is not necessary, as it was proved in \cite{PTY}.

So, now we deal with the Mellin transform of residue function (\ref{f1}) which is given by the integral
\begin{equation}\label{f3}
s\mapsto \Gamma_{f}(s,\varphi) =\int\limits_{{\mathbb R}_{+}^{p}}I(\varepsilon ,\varphi)\varepsilon ^{s-I}d\varepsilon,
\end{equation}   
where $s=(s_1,\dots , s_p)\in {\mathbb C}^p$ and $\varepsilon ^{s-I}d\varepsilon :=\varepsilon_{1}^{s_1-1}\cdot \ldots \cdot \varepsilon_{p}^{s_p -1}d\varepsilon_1\wedge \ldots \wedge d\varepsilon_p .$ It admits a representation as the following solid volume integral
\begin{equation}\label{f4}
\Gamma_{f}(s,\varphi)=\frac{c_p}{(p-1)!}\int\limits_{V}^{}\prod\limits_{k=1}^{p}|f_k|^{2(s_k-1)}\bigwedge\limits_{k=1}^{p}\overline{d f_k}\wedge\varphi .
\end{equation}
Theorem 2.2 in \cite{PTY} states that the current  $T_{f}(\varphi)$ can be given in terms of the function $\Gamma_{f}(s,\varphi)$ as follows
\begin{equation}\label{f5}
T_{f}(\varphi)=\lim\limits_{\tau \to 0+} \frac{1}{(2\pi i)^p}\int\limits_{\gamma+i{\mathbb R}^p} \tau^{-|s|}\Gamma (|s|+1)\prod\limits_{k=1}^{p}\Gamma (1-s_k) \Gamma_{f}(s, \varphi) ds ,
\end{equation}
where  $\gamma\in (0;1)^p$ and $|s|:=s_1+\ldots +s_p$. 

According to \cite[Th. 2]{PT95}, the  Mellin transform $\Gamma_{f}(s,\varphi)$ is holomorphic on the set
$\{s\in {\mathbb C}^{p}: {\Re s}\in {\mathbb R}_{+}^{p} \}$ and it has a meromorphic continuation to all of ${\mathbb C}^p$.
In particular, for a monomial mapping
\begin{equation}\label{f4.1}
\zeta^{ A}=(\zeta^{\alpha_1},\dots , \zeta^{\alpha_p}):{\mathbb C}^n\to {\mathbb C}^p,
\end{equation}
given by a $p\times n$ -- matrix $A$, whose row vectors are $\alpha_j\in {\mathbb Z}_{\geqslant}^{n}$, the Mellin transform $\Gamma_{\zeta^A}(s,\varphi)$ has simple poles on hypersurfaces $\langle \alpha^k, s \rangle=-\nu$, $\nu\in {\mathbb N}\cup \{0\}$, $1\leqslant k\leqslant n$, where $\alpha^k$ denote column vectors of the matrix $A$. Moreover, near the origin one has the following representation
\begin{equation}\label{f7} 
\Gamma_{\zeta^A}(s,\varphi)=\sum\limits_{\sharp K=p}^{} \frac{c_{K}}{\langle \alpha^{k_1}, s \rangle \cdot \ldots \cdot \langle \alpha^{k_p}, s\rangle } + Q(s)
\end{equation}
with constants $c_K$, where $K=(k_1,\dots , k_p)$ is an ordered set of indices, and $Q(s)$ is a finite sum of functions with simple poles along fewer than $p$ hyperplanes. 

First of all, we clarify the view of integral (\ref{f4}) in the monomial case.

\begin{lm}
Let $\zeta^{A}$  be a monomial mapping.  The Mellin transform
$\Gamma_{\zeta^{A}}(s,\varphi)$ is equal to the sum
$$
\sum\limits_{\sharp I=p}^{}\Gamma_{\zeta^{A}}^{I}(s,\varphi)
$$
over all ordered sets  $I=(i_1,\ldots , i_p)$, $1\leqslant i_1 <\ldots <i_p\leqslant n$, whose summands $\Gamma_{\zeta^{A}}^{I}(s,\varphi)$ are defined as follows
\begin{equation}\label{f4.2}
\begin{split}
\frac{c_I{\Delta_I}}{(2\pi i)^p}\int\limits_{{\mathbb C}^n}^{}\varphi_{I} (\zeta) \bigwedge\limits_{j\in I}\left(\frac{|\zeta_j|^{2(\langle \alpha^{j}, s\rangle -1)}}{\zeta_{j}^{|\alpha^{j}|-1}}d\bar{\zeta}_j\wedge d\zeta_{j}\right)
\bigwedge\limits_{k\notin I}\left(\frac{|\zeta_k|^{2\langle \alpha^{k}, s\rangle}}{\zeta_{k}^{|\alpha^{k}|}}
 d\bar{\zeta}_k\wedge d\zeta_{k}\right).
\end{split}
\end{equation}
In (\ref{f4.2}) $\Delta_{I}$ denotes a minor of the matrix $A$, $|\alpha^{j}|=\alpha_1^{j}+\ldots +\alpha_p^{j}$ is the sum of $j^{,}$th column vector entries, $\varphi_{I}(\zeta)$ is a coefficient of the test form $\varphi=\sum\limits_{I}\varphi_{I}d\zeta \wedge d \bar{\zeta} [I]$ and $c_I=(-1)^{\vert I \vert - (p+1)p/2+(n-p+1)(n-p)/2}$.
\end{lm}

{\it Proof.} Let us remark that
\begin{equation}\label{f4.3}
\prod\limits_{k=1}^{p}|\zeta^{\alpha_k}|^{2(s_k-1)}=\prod\limits_{k=1}^{n}\left(|\zeta_k|^2\right)^{\langle \alpha^k , s\rangle -|\alpha^{k}|},\end{equation}
\begin{equation}\label{f4.4}
\bigwedge\limits_{k=1}^{p}\overline{d\zeta^{\alpha_k}}=\bar{\zeta}^{\alpha_{1}}\cdot\ldots\cdot \bar{\zeta}^{\alpha_p}\sum\limits_{\sharp I=p}^{}\Delta_I \frac{d\bar{\zeta}_{I}}{\bar{\zeta}_{I}},
\end{equation}
where $\bar{\zeta}_I=\bar{\zeta}_{i_1}\cdot \ldots \cdot \bar{\zeta}_{i_p}$, $d\bar{\zeta}_I=d\bar{\zeta}_{i_1}\wedge \ldots \wedge d\bar{\zeta}_{i_p}$ and $\Delta_I=\mbox{{\rm det}}A_{I}$, $A_{I}=(\alpha^{j})$, $j\in I$. We insert expressions (\ref{f4.3}), (\ref{f4.4}) in to integral (\ref{f4}), whose elementary transformations complete the proof.


Lemma 1 and formula (\ref{f5}) imply that the current $T_{\zeta^A}$  is equal to the sum
$$
T_{\zeta^A}(\varphi)=\sum\limits_{\sharp I=p}^{}T_{\zeta^A}^{I}(\varphi),
$$
whose summands can be represented as follows:

\begin{equation}\label{f9}
T_{\zeta^A}^{I}(\varphi)=\lim\limits_{\tau \to 0+} \frac{1}{(2\pi i)^p}\int\limits_{\gamma+i{\mathbb R}^p}\tau^{-|s|}\Gamma (|s|+1)\prod\limits_{k=1}^{p}\Gamma (1-s_k) \Gamma_{\zeta^A}^I(s, \varphi) ds. 
\end{equation}
In fact, the set of admissible values for $\gamma$ in (\ref{f9}) may be enlarged to the set
$$
M=\left(\mbox{\it{int}} K^{I}\right)\cap \left\{x\in {\mathbb R}^p: x_1<1,\dots , x_p<1\right\},
$$
where $\mbox{\it {int}} K^{I}$ denotes the interior of the cone $K^{I}$ whose definition is given in (\ref{f7.1}) below.

\section{Structure Theorem}

We consider the Bochner--Martinelli  residue current, corresponding to the system of monomials (\ref{f4.1}). Let us remind that $\alpha^j$ is the $j$'th vector column of the matrix $A$ and   $|\alpha^j|$ denotes the sum of its entries. From Lemma 1 it follows that the Mellin transform    $\Gamma_{\zeta^A}(s, \varphi)$ is equal to the sum of functions $\Gamma_{\zeta^A}^{I}(s, \varphi)$, each of which is represented by integral (\ref{f4.2}). Let us fix a set $I=(i_1,\dots , i_p)$.  The corresponding function  $\Gamma_{\zeta^A}^{I}(s, \varphi)$ has poles on $p$ families of hyperplanes
\begin{equation}\label{f7}
\langle\alpha^{k}, s\rangle =-\nu ,\,\, \nu \in {\mathbb N}\cup \{0\},\,\, k\in I,
\end{equation}
(see \cite{PT95}).
In the space ${\mathbb R}^p$ with variables  $x_j={\Re}s_j$ we introduce the cone
\begin{equation}\label{f7.1}
K^{I}=\{x\in{\mathbb R}^p: \left\langle \alpha^{k} , x\right\rangle\geqslant 0,
\, k\in I\}
\end{equation}
and the intersection
$$
K_0^{I}=K^{I}\cap \Pi_{\leqslant}
$$
of this cone with the half-space
$$
\Pi_{\leqslant}:=\{x\in{\mathbb R}^p: x_1+\dots +x_p \leqslant 0\}.
$$
Denote by $q=q(I)$ the codimension of the cone $K_0^{I}$. It follows that $K_{0}^{I}$ admits the representation
$$
K_0^{I}=\{x\in {\mathbb R}^p:\left\langle \alpha^{j} , x\right\rangle = 0,\, j\in J, \left\langle \alpha^{k} , x\right\rangle\geqslant 0,\,\, k\in I\setminus J\},
$$
where $J=J(I)$ is a subset of $I$ the cardinality $q$. If $q=0$, then the set $K_0^I$ contains interior points of the half--space $\Pi_{\leqslant}$. If $q\geqslant 1$, then $K_0^I$ is a $(p-q)$--dimensional  face of the cone  $K^I$, lying on the boundary of the set $\Pi_{\leqslant}$.

\begin{theo}
The Bochner-Martinelli residue current $T_{\zeta^A}$, associated with the mapping $\zeta^A$, admits the representation
$$
T_{\zeta^A}=\sum\limits_{\sharp I=p}^{}T_{\zeta^A}^{I},
$$
where
\begin{equation}\label{f8}
T_{\zeta^A}^{I}=c_I \bigwedge\limits_{j\in J}^{}\bar\partial \left[\frac{1}{\zeta_{j}^{|\alpha^{j}|}}\right]
\bigwedge_{k\in I\setminus J}^{}\left(\frac{1}{\zeta_k^{|\alpha^k|}}\frac{d\bar\zeta_k}{\bar\zeta_k}\right)\bigwedge_{l\notin I}^{}
\left[\frac{1}{\zeta_{l}^{|\alpha^{l}|}}\right]
F(|\zeta_1|^2,..[J]..,|\zeta_n|^2).
\end{equation}
The $F(|\zeta_1|^2,..[J]..,|\zeta_n|^2)$ is a hypergeometric function (whose Mellin--Barnes integral representation is given in formula (\ref{f12}) below) and the constant $c_I$ is defined in Lemma 1. In particular, if $J=I$, then 
\begin{equation*}
T_{\zeta^A}^{I}=c_I \bigwedge\limits_{j\in I}^{}\bar\partial \left[\frac{1}{\zeta_{j}^{|\alpha^{j}|}}\right]
\bigwedge_{l\notin I}^{}
\left[\frac{1}{\zeta_{l}^{|\alpha^{l}|}}\right].
\end{equation*}

\end{theo}

\begin{rem}
The support of the current $T_{\zeta^A}^{I}$ has complex codimension $q(I)$.
\end{rem}

\section{Leray Residue Form}

To prove representation (\ref{f8}) we need to decrease the number of integrations and take the limit in (\ref{f9}). It is achieved by calculations of the composite residue form of the integrand in (\ref{f9}). The existence of the residue form  for a given closed regular form was proved by Leray. Let us remind the definition of it (see \cite[Ch. III]{AYu}). Let $L\subset {\mathbb C}^n$ be a complex analytic submanifold of dimension $n-1$ given in some neighborhood of a point $a\in L$ by the equation $\lambda (z)=0$, where $\lambda (z)$ is a holomorphic function. Suppose that $\psi (z)$ is a closed regular form on ${\mathbb C}^n\setminus L$ that has on $L$ a polar singularity of the first order. It means that the form $\psi (z)\lambda (z)$ extends to some neighborhood of $a\in L$ as a $C^{\infty}$--form. Then in some neighborhood $U_a$ of $a\in L$ there exist regular forms $\chi (z)$, $\theta (z)$ such that in $U_a$
\begin{equation}\label{f8.1}
\psi (z)=\frac{d\lambda (z)}{\lambda (z)}\wedge \chi (z) +\theta (z).
\end{equation}
The closed form on $L$ defined in a neighborhood of each point $a\in L$ by the restriction $\chi (z)\vert_{L}$ is called the residue form of the form $\psi$ and denoted by ${\mbox{\rm res}}[\psi]$. 

For instance, if $L=\{z\in {\mathbb C}^n : \lambda (z)=0\}$ and
\begin{equation*}
\psi (z)=\frac{f(z)}{\lambda (z)}dz_1\wedge\ldots\wedge dz_n,
\end{equation*} 
where functions $f(z)$, $\lambda(z)$ are holomorphic in ${\mathbb{C}^n}$ and ${\mbox{\rm grad}}\lambda (z)\vert_{L}\ne 0$, then
$$
{\mbox{\rm res}}[\psi]=(-1)^{j-1}\left.\frac{f(z)}{\lambda_{z_j}^{'}(z)}\right\vert_{L}dz_1\wedge ..[j]..\wedge dz_n 
$$
in points where $\lambda_{z_j}^{'}(z)\ne 0$.

Let us consider complex analytic $(n-1)$--dimensional submanifolds $L_1,\dots , L_q\subset {\mathbb C}^n$ which intersect in common points transversally. It means that in some neighborhood of a common point $a$ holomorphic functions defining locally submanifolds $L_1$,..., $L_q$ have linearly independent gradients.   
If the closed form $\psi (z)$ has simple poles on submanifolds $L_1, \ldots , L_q$, then repeated application of (\ref{f8.1}) enables us to define the form
$$
{\mbox{\rm res}}^{q}[\psi]={\mbox{\rm res}}_{q}\circ \ldots \circ {\mbox{\rm res}}_{1}[\psi],
$$
called the composite residue form.

We remark that if
$$
\psi=\frac{d\lambda_1\wedge \ldots \wedge d\lambda_q \wedge \chi}{\lambda_1\cdot\ldots\cdot \lambda_q}
$$
then 
$$
{\mbox{\rm res}}^{q}[\psi]=\chi\vert_{L_1\cap\ldots\cap L_q},
$$
where $\lambda_j=\lambda_j(z)$ are holomorphic functions defining locally submanifolds $L_j$, $j=1,\dots , q$, ${\mbox{\rm grad}}\lambda_{j}\vert_{L_1\cap\ldots\cap L_q}\ne 0$.

Let us consider integral (\ref{f9}) whose integrand we denote by $\omega_{\tau}^{I}(s)$. If we choose $\gamma$ in the set $\Pi_{\leqslant}$, i.e. $|\gamma|=\gamma_1+\ldots +\gamma_p<0$, then
$$
\int\limits_{\gamma+i{\mathbb R}^p}^{}\omega_{\tau}^{I}(s) \to 0,\,\, \tau \to 0+,
$$
since the form $\omega_{\tau}^{I}(s)$ contains the factor $\tau^{-|s|}$.

Having in mind this fact we can distinguish  polar hyperplanes 
$$
L_j=\{s\in{\mathbb C}^p: \langle\alpha^j , s\rangle =0\},\,\, j\in J=(j_1,\ldots , j_q),
$$
of the form $\omega_{\tau}^{I}(s)$ with respect to whose we are able to compute the residue form
$$
\mbox{\rm{res}}_{L^q}\left[\omega_{\tau}^{I}(s)\right],
$$
where $L^{q}:=\bigcap\limits_{j\in J}^{}L_j$, $\Re (L^{q})\supset K_{0}^{I}$.

Let $q=1$, then $L^{1}=L_{j_1}$. We consider the ray
$$
l=\gamma +\{\langle \alpha^{j_1}, x\rangle\leqslant 0, \langle\alpha^{j}, x\rangle=0,\, j\in I, j\ne j_1\},
$$
intersecting the set $\Re L^1$ in the point $\gamma^1$. For any points $\gamma$ and $\varkappa$ of $l$ lying on different sides of $\gamma^1$, the Cauchy formula yields the following equality
\begin{equation}\label{f8.2}
\int\limits_{\gamma+i{\mathbb R}^p}^{}\omega_{\tau}^{I}(s)=\int\limits_{\varkappa+i{\mathbb R}^p}^{}\omega_{\tau}^{I}(s)+2\pi i
\int\limits_{\gamma^1+i\Im L^1}^{} \mbox{\rm {res}}_{L^1}\left[\omega_{\tau}^{I}(s)\right].
\end{equation}
If $\varkappa\in\Pi_{\leqslant}$, then the first integral in (\ref{f8.2}) tends to zero as $\tau\to 0+$. The second one in the right hand side has the same structure as the integral in the left-hand side, but in subspace $L_{j_1}$.

Having $q>1$, we repeat $q-1$ times again residue theorem. As a result, we obtain the identity
$$
\int\limits_{\gamma+i{\mathbb R}^p}^{}\omega_{\tau}^{I}(s)=(2\pi i)^q \lim\limits_{\tau\to 0+}^{}\int\limits_{\gamma^q+i\Im L^q}^{}{\mbox{\rm res}}_{L^q}\left[\omega_{\tau}^{I}(s)\right]+o(\tau),
$$
where $\gamma^q\in\Re L^q$. The set $L^q\subset \{s:|s|=0\}$, hence 
$$
\mbox{\rm res}_{L^q}\left[\omega_{\tau}^{I}(s)\right]=\prod\limits_{j=1}^{p}
\Gamma (1-s_j){\rm res}_{L^q}\left[\Gamma_{\zeta^A}^{I}(s,\varphi)ds\right].
$$
Therefore, we have the following statement which was proved in \cite{PTY} for $p=n$.

\begin{lm}
The limit of the $p$--fold integral (\ref{f9}) may be written as the $(p-q)$--fold integral
\begin{equation}\label{f10}
T_{\zeta^A}^{I}(\varphi)=\frac{1}{(2\pi i)^{p-q}}\int\limits_{\gamma^{q}+i{\rm \Im} L^{q}}^{}\prod\limits_{j=1}^{p}
\Gamma (1-s_j){\rm res}_{L^q}\left[\Gamma_{\zeta^A}^{I}(s,\varphi)ds\right],
\end{equation}
where $\gamma^{q}\in M^{q}:=\mbox{\it{rel int}}(K_{0}^{I}\cap \overline{M})$.
\end{lm}
So, we need to compute the residue of the form $\Gamma_{\zeta^A}^{I} (s, \varphi)ds$, where the function  $\Gamma_{\zeta^A}^{I} (s, \varphi)$ is given by integral (\ref{f4.2}). For that let us introduce  new coordinates $\lambda=(\lambda_j)=sA$, $\lambda_j=\langle\alpha^j , s\rangle$, $j=1,\ldots , n$. Letting $A_I$ denote the square matrix, whose column vectors are $\alpha^j$, $j\in I$, we thus get
$$
sA_I=\lambda_I,\,\, s=\lambda_I A_{I}^{-1},
$$
where $\lambda_I=(\lambda_{i_1},\ldots , \lambda_{i_p})$, $\Delta_I={\mbox {\rm det} A_I}\ne 0$. If $\Delta_I=0$, then $\Gamma_{\zeta^A}^{I} (s, \varphi)=0$. Integral (\ref{f4.2}) in coordinates $\lambda_I$ looks as follows:
\begin{equation}\label{f11}
\frac{c_I\Delta_I}{(2\pi i)^p} 
\int\limits_{{\mathbb C}^n}^{}\varphi_{I} (\zeta) \bigwedge\limits_{j\in I}\left(\frac{|\zeta_j|^{2(\lambda_j -1)}}{\zeta_{j}^{|\alpha^{j}|-1}}d\bar{\zeta}_j\wedge d\zeta_{j}\right)
\bigwedge\limits_{k\notin I}\left(\frac{|\zeta_k|^{2\langle \lambda_I, \mu^{k} \rangle}}{\zeta_{k}^{|\alpha^{k}|}} 
d\bar{\zeta}_k\wedge d\zeta_{k}\right),
\end{equation}
where
$$
\mu^{k}=\sum\limits_{l=1}^{p}\alpha_{l}^{k}\beta^l
$$
and $\beta^{l}$ is the $l$'th column vector of the inverse matrix $A_{I}^{-1}$. 

Now let perform the integration in (\ref{f11}) with respect to variables $\zeta_J=(\zeta_{j_1},\dots , \zeta_{j_q})$ using polar coordinates. To this end for the smooth function $\varphi_{I}(\zeta):{\mathbb C}^n\to {\mathbb C}$ and $\alpha=(|\alpha^{j_1}|-1,\ldots , |\alpha^{j_q}|-1)\in {\mathbb N}^q$ we consider a decomposition of the type
\begin{equation}\label{f11.1}
\varphi_{I}(\zeta)=\sum\limits_{j\in J}^{}\sum\limits_{k+l<|\alpha^{j}|-1}^{}\psi_{kl}^{j}(\zeta [j])\zeta_{j}^{k}\bar{\zeta}_{j}^{l}+
\sum\limits_{K+L=\alpha}^{}\psi_{KL}(\zeta)\zeta_{J}^{K}\bar{\zeta}_{J}^{L},
\end{equation}
where all the coefficient functions are smooth and the $\psi_{kl}^{j}$ are independent of the variable $\zeta_j$, see \cite{PT95}. Decomposition (\ref{f11.1}) is a result of the Taylor expansion for the function $\varphi_{I}(\zeta)$ at the origin with respect to variables $\zeta_J=(\zeta_{j_1},\dots , \zeta_{j_q})$. Introducing polar coordinates $\zeta_j=r_j e^{i\theta_j}$, $j\in J$, in (\ref{f11}) and using (\ref{f11.1}) we get two series of inner integrals
\begin{equation}\label{f11.2}
\int\limits_{0}^{2\pi} e^{i(k-l+1-|\alpha^{j}|)\theta_j}d\theta_j,
\end{equation}
and
\begin{equation}\label{f11.3}
\int\limits_{0}^{2\pi}e^{i(K_j-L_j+1-|\alpha^{j}|)\theta_j}d\theta_j.
\end{equation}
Integrals (\ref{f11.2}) are all equal to zero, since $k-l\geqslant k+l < |\alpha^{j}|-1$. In series (\ref{f11.3}) the only one integral is nonzero when $L_j=0$, $K_j=|\alpha^{j}|-1$. Moreover, remind that the test form $\varphi$ has compact support, so, in fact, we integrate with respect to a variable $r_j$ over a segment $[0;R_j]$.

As a result of integration in (\ref{f11}) with respect to variables $\zeta_J=(\zeta_{j_1},\dots , \zeta_{j_q})$, we have the following representation
\begin{equation}\label{f12}
\begin{split}
\Gamma_{\zeta^A}^{I}=\frac{c_I\Delta_I}{(2\pi i)^{p-q}}\frac{1}{\lambda_{j_1}\ldots \lambda_{j_q}} 
\left\{\int\limits_{{\mathbb C}^{n-q}} \varphi_{I}^{\alpha}(\zeta [J])\bigwedge \limits_{j\in I\setminus J}^{} \left(
\frac{|\zeta_j|^{2\lambda_j}}{\zeta_j^{|\alpha^{j}|}}\frac{d\bar{\zeta_j}}{\bar{\zeta_j}}\wedge d\zeta_j\right) \right.\\
\left.\bigwedge\limits_{k\notin I}^{}\left(\frac{|\zeta_k|^{2\langle\lambda_I, \mu^{k} \rangle}}{\zeta_k^{|\alpha^{k}|}}d\bar{\zeta_k}\wedge d\zeta_k\right) 
+{\tilde Q}(\lambda_I)\right\},
\end{split}
\end{equation}
where
$$
\varphi_{I}^{\alpha}(\zeta [J])=\frac{1}{(|\alpha^{j_1}|-1)! \dots (|\alpha^{j_q}|-1)!}\frac{\partial^{|\alpha^{j_1}|+\dots +|\alpha^{j_q}|-q}}
{\partial\zeta_{1}^{|\alpha^{j_1}|-1} \dots \partial\zeta_{q}^{|\alpha^{j_q}|-1}}\varphi_{I}(0_{J},\zeta [J]),
$$
$\zeta[J]=(\zeta_1,..[J].., \zeta_n)$ and ${\tilde Q}(\lambda_I)$  is a holomorphic function in a neighborhood of the origin, belonging to the ideal $\left\langle \lambda_{j_1},\dots , \lambda_{j_q}\right\rangle$.

The form $\Gamma_{\zeta^A}^{I} ds$ in new coordinates $\lambda_I$ may be written as follows:
\begin{equation}\label{f13}
\Gamma_{\zeta^A}^{I} ds=\bigwedge\limits_{j\in J}^{}\frac{d\lambda_j}{\lambda_j} \wedge \left(R(\lambda_I)d \lambda_{I}[J]\right)
+Q(\lambda_I),
\end{equation}
where $\lambda_{I}[J]=(\lambda_{i_1},.. [J] .., \lambda_{i_p})$, $d\lambda_{I}[J]=d\lambda_{i_1} \wedge .. [J] .. \wedge d\lambda_{i_p}$, $Q(\lambda_I)$ is a meromorphic $p$--form,
which has poles on  $r < q$ coordinate hyperplanes. The restriction of the form $R(\lambda_I)d\lambda_{I}[J]$ on the set $L^q$ gives the residue form that we need:
\begin{equation}\label{f13.1}
{\mbox{\rm res}}_{L^q}\left[\Gamma_{\zeta^A}^{I} ds\right]= 
R(0_J,\lambda_{I}[J])d \lambda_{I}[J],
\end{equation}
where
\begin{eqnarray*}
R(0_J,\lambda_{I}[J])=\frac{c_I}{(2\pi i)^{p-q}}\int\limits_{{\mathbb C}^{n-q}} \varphi_{I}^{\alpha}(\zeta [J])
\left(\bigwedge \limits_{j\in I\setminus J}^{}
\frac{|\zeta_j|^{2\lambda_j}}{\zeta_j^{|\alpha^{j}|}}\frac{d\bar{\zeta_j}}{\bar{\zeta_j}}\wedge d\zeta_j\right)\bigwedge \\
\left(\bigwedge\limits_{k\notin I}^{}\frac{|\zeta_k|^{2\langle\lambda_{I}[J], \mu^{k}[J]\rangle}}{\zeta_k^{|\alpha^{k}|}}d\bar{\zeta_k}\wedge d\zeta_k\right),
\end{eqnarray*}
and  $\mu^{k}[J]$ means that coordinates $\mu_{j}^{k}$, $j\in J$, in $\mu^{k}$ are removed.

\section{Proof of Theorem 1}
We have now collected almost all elements for the proof of Theorem 1.
If $q=0$, then  $T_{\zeta^A}^{I}=0$. In this case  ${\rm dim}_{{\mathbb R}} K_{0}^I =p$ and we can choose $\gamma$ in the set $M$ so that
$|\gamma|=\gamma_1 +\dots +\gamma_p <0$. Consequently, the restriction of the integrand in (\ref{f9}) to ${\gamma +i{\mathbb R}^p}$ tends to zero as $\tau \to 0+$, and hence, $T_{\zeta^A}^{I}=0$.

Assuming that $q\geqslant 1$, we apply Lemma 2.  The set $L^q$ admits the following parametrization:
$$
s_{l}=\sum\limits_{j\in I\setminus J}^{}\beta_{j}^{l}\lambda_{j},\,\, l=1,\dots , p.
$$
Substituting the residue form (\ref{f13.1}) into formula (\ref{f10}) and applying Fubini's theorem, we can observe the structure of the residue current $T_{\zeta^A}^{I} (\varphi)$. It admits the following representation:
\begin{eqnarray*}
\frac{c_I}{(2\pi i)^{p-q}}
\int\limits_{{\mathbb C}^{n-q}} \varphi_{I}^{\alpha}(\zeta [J])
\bigwedge\limits_{j\in I\setminus J}^{}\left(
\frac{1}{\zeta_j^{|\alpha^{j}|}}\frac{d\bar{\zeta_j}}{\bar{\zeta_j}}\wedge d\zeta_j\right) \times \nonumber \\
\bigwedge\limits_{k\notin I}^{}\left(\frac{1}{\zeta_k^{|\alpha^{k}|}}d\bar{\zeta_k}\wedge d\zeta_k\right) F(|\zeta_1|^2,..[J]..,|\zeta_n|^2),
\end{eqnarray*}
where the function $F(|\zeta_1|^2,..[J]..,|\zeta_n|^2)$ is representable as the Mellin--Barnes integral
\begin{equation}\label{f12}
\frac{1}{(2\pi i)^{p-q}}\int\limits_{i{\mathbb R}^{p-q}}^{}\prod\limits_{l=1}^{p}\Gamma \left(1-\sum\limits_{j\in I\setminus J}^{}\beta_{j}^{l}\lambda_{j}\right)
\prod\limits_{j\in I\setminus J}^{}\left(|\zeta_j|^2\prod\limits_{k \notin I}^{}|\zeta_k|^{2\left\langle \alpha^k, \beta_j\right\rangle}\right)^{\lambda_j}d\lambda_{I}[J],
\end{equation}
where $\beta_j$ is the $j$'th row of the matrix $A_I^{-1}$. 

If $J=I$, then the integration set in (\ref{f10}) is zero-dimensional, so the current $T_{\zeta^{A}}^{I}(\varphi)$ admits the following representation
\begin{eqnarray*}
c_{I}
\int\limits_{{\mathbb C}^{n-p}} \varphi_{I}^{\alpha}(\zeta [I])
\bigwedge\limits_{k\notin I}^{}\left(\frac{1}{\zeta_k^{|\alpha^{k}|}}d\bar{\zeta_k}\wedge d\zeta_k\right)=
c_I \bigwedge\limits_{j\in I}^{}\bar\partial \left[\frac{1}{\zeta_{j}^{|\alpha^{j}|}}\right]
\bigwedge_{k\notin I}^{}
\left[\frac{1}{\zeta_{k}^{|\alpha^{k}|}}\right]. 
\end{eqnarray*}
This concludes the proof of Theorem 1.

}
\end{document}